\documentclass[12pt]{amsart}
 \usepackage{amsmath,amsthm}

\def\P{{\mathbb P}}
\def\Z{{\mathbb Z}}
\def\C{{\mathbb C}}
\def\N{{\mathbb N}}
\def\Q{{\mathbb Q}}

\def\I{{\mathcal J}}
\def\OOO{{\mathcal O}}

\def\T{\Theta}
\def\e{\varepsilon}

\theoremstyle{plain}
\newtheorem{thm}{Theorem}[section]

\newtheorem{prop}[thm]{Proposition}
\newtheorem{cor}[thm]{Corollary}

\newtheorem{claim}[thm]{Claim}

\theoremstyle{definition}
\newtheorem{df}[thm]{Definition}
\newtheorem{rem}[thm]{Remark}
\newtheorem{ex}[thm]{Example}

\title{Multiplier Ideals in Algebraic Geometry}
\author{Samuel Grushevsky}
\address{Mathematics Department, Princeton University, Fine Hall,
Washington Road, Princeton, NJ 08544, USA}
\email{sam@math.princeton.edu}
\date{December 19, 2004}
\thanks{Partially supported by the NSF Mathematical Sciences
Postdoctoral Research Fellowship} 

\begin{document}
\maketitle

\section*{Introduction}
In this introductory survey text we introduce multiplier ideal
sheaves in the context of general vanishing theorems and
log-resolution of singularities. After discussing some basic
properties of multiplier ideals, we then follow \cite{lazmult} to
obtain Koll\'ar's bound from \cite{kolmult} on the multiplicity of
theta divisors on abelian varieties. We then develop the theory of
asymptotic multiplier ideals and expose some of the ideas involved
in the algebraic interpretation of Siu's proof of deformation
invariance of plurigenera for varieties of general type
\cite{siuplur}. We also define Nadel multiplier ideals in the
analytic setting, and explain some ideas behind the analytic
proofs of the vanishing theorems and of the deformation invariance
of plurigenera, to give the reader an idea of the analytic side of
the theory.

This text is by no means meant to be a complete self-contained
introduction to the broad and rich field of multiplier ideals, and
the reader is encouraged to look at the excellent rigorous in-depth
treatments of the subject like the ones in \cite{lazbook},
\cite{siunote}, \cite{dembook} and references therein. Thus the
purpose of this article is to explain some of the ideas and
techniques in multiplier ideals, and to encourage the reader to
learn this exciting field in more detail. Our exposition mostly
follows \cite{lazbook} for the algebraic and \cite{siunote} for the
analytic story.

{\bf Acknowledgements.} I owe a debt of gratitude to Yum-Tong Siu,
from whose classes (including the one the notes for which constitute
\cite{siunote}), talks, and most importantly many insightful
conversations with whom I first learned the subject. I have greatly
benefited from detailed discussions with Robert Lazarsfeld, and
learned the algebraic side of the theory from the preliminary
versions of his book \cite{lazbook}. Professor Lazarsfeld has also
made numerous valuable comments and suggested many clarifications
for a preliminary version of this text. I am grateful to Carolina
Araujo and Jordan Ellenberg, who attended and commented on a
preliminary version of the talk, which led to further improvements
in the talk and in this paper. I would also like to thank Richard
Thomas for reading a draft of this text very closely and making many
useful suggestions.

\section{Classical theory and Kodaira's vanishing}
{\bf Convention.} {\it We will work over the field of complex
numbers, and for simplicity will assume all the varieties to be
smooth, though most of the methods have been generalized to deal
with the singular case as well.}

In this text we will be concerned with the study of line bundles on
a projective variety $X$ of (complex) dimension $n$ and their
cohomology. We do not make any distinction between a line bundle and
the corresponding divisor, and use additive notations for line
bundles, i.e. denote $L\otimes L$ by $2L$. We denote the space of
sections of a bundle $L$ over $X$ by $\Gamma(X,\OOO_X(L))$; its
dimension is denoted $h^0(X,\OOO_X(L))$. The complete linear system
$|L|$ is the space of all one-dimensional linear subspaces of
$\Gamma(X,\OOO_X(L))$ --- the elements of $|L|$ are divisors
linearly equivalent to $L$; the canonical bundle of $X$ is denoted
by $K_X$. Let us now start with some classical definitions and
results.

\begin{df} A line bundle $L$ on $X$ is called very ample if its
sections embed $X$ into a projective space, i.e. if for any basis
$s_0\ldots s_N\in \Gamma(X,\OOO_X(L))$ the map $X\to\P^N$ obtained
by sending a point $x\in X$ to $(s_0(x):s_1(x):\ldots:s_N(x))$ is a
well-defined embedding. A line bundle is called ample if there
exists some number $m\in\N$ such that $mL$ is very ample.
\end{df}

The following numerical criterion of ampleness is classical:

\begin{thm}[Nakai-Moishezon-Kleiman criterion] A line bundle $L$ on
$X$ is ample if and only if for any subvariety $Y\subset X$ of any
dimension $d$ the intersection $Y\cdot L^d>0$.
\end{thm}

Ample line bundles are very special from the point of view of the
cohomology theory in view of the following

\begin{thm}[Kodaira's vanishing]
If $L$ is ample, then
$$
 H^i(X,\OOO_X(K_X+L))=H^{n-i}(X,\OOO_X(-L))=0
$$
for all $i>0$ (the two cohomology groups are dual by Serre
duality).
\end{thm}

One may wonder whether this cohomological vanishing property
characterizes ample line bundles. The answer is no, and the right
question to ask instead is for which bundles the vanishing holds. It
seems the vanishing should hold at least for limits of ample
bundles, which may not necessarily be ample themselves. These limits
are well understood, and they are defined by relaxing both
conditions of ampleness.

\begin{df} A line bundle $L$ is called numerically effective, or nef, if
for any curve $C\subset X$ one has $L\cdot C\ge 0$. A line bundle is
called big if for some $m\in\Z$ the rational map from $X$ to the
projective space $\P^{h^0(X,\OOO_X(mL))-1}$ given by sections is
birational. A variety $X$ is said to be of general type if its
canonical bundle $K_X$ is big.
\end{df}

The reason one uses only intersections with curves, and not with
higher-dimensional subvarieties, in the definition of nefness is the
following

\begin{thm}[Kleiman's theorem] If a line bundle $L$ is nef, then its
degree on any $d$-dimensional subvariety $Y\subset X$ is
non-negative, i.e. $L^d\cdot Y\ge 0$.
\end{thm}

This theorem means that for testing whether a divisor is nef
intersecting only with curves suffices. One may thus wonder whether
for testing ampleness checking the fact that $L\cdot C>0$ for all
curves suffices, and the answer is no, with the counterexamples
already present for some smooth projective surfaces --- see
\cite{lazbook}, 1.2-1.5 for a more detailed discussion.

\begin{thm}[Kawamata-Viehweg]
The vanishing theorem holds for nef and big bundles, i.e. if $L$ is
both big and nef, then $H^i(X,\OOO_X(K_X+L))=0$ for $i>0$.
\end{thm}

Kodaira's vanishing theorem has found innumerable applications and
has been used to obtain results in a variety of settings in
algebraic geometry. The culmination of our discussion here will be
an outline of the ideas behind the proof of the deformation
invariance of plurigenera, which we now state.

\begin{thm}[Siu, \cite{siuplur}] The plurigenera of a variety of general type
are deformation invariant, i.e. if we have a family $X\to \Delta$ of
$n$-dimensional varieties over the unit complex disk $\Delta$, with
all fibers $X_t$ being of general type, then for any integer $m$ the
$m^{\text th}$ plurigenus $h^0(X_t,\OOO_{X_t}(mK_{X_t}))$ is
independent of $t$ for $t$ small enough.
\end{thm}

It is quite easy to see that plurigenera are upper semicontinuous,
because any family of sections of $mK_{X_t}$ that exists for all
$t\ne 0$ can be extended to give an element of
$\Gamma(X_0,\OOO_{X_0}(mK_{X_0})$. Thus the hard part of the proof
is to show that the plurigenus cannot accidentally increase for
$X_0$, i.e. essentially that any section of $mK_{X_0}$ gives rise to
a section in $\Gamma(X_t,\OOO_{X_t}(mK_{X_t}))$ for all $t$
sufficiently small. By gluing these together for all $t$, this is
the same as asking whether an element of
$\Gamma(X_0,\OOO_{X_0}(mK_{X_0}))$ extends to a section in
$\Gamma(X,\OOO_X(mK_X))$.

At first sight it does not seem that the vanishing theorems are
related to the invariance of plurigenera. However, vanishing
theorems can in fact yield the invariance of plurigenera directly in
some cases:

\begin{prop}\label{ampleplurigenus}
If $K_{X_0}$ is big and nef, the invariance of plurigenera holds.
\end{prop}
\begin{proof}
Indeed, if $K_{X_0}$ is big and nef, all its multiples $mK_{X_0}$
are also big and nef, and so all the higher cohomologies
$H^i(X_0,\OOO_{X_0}(mK_{X_0}))$ are zero. By semicontinuity it
follows that $H^i(X_t,\OOO_{X_t}(mK_{X_t}))=0$ for all $t$ small
enough. Thus we see that $h^0(X_t,\OOO_{X_t}(mK_{X_t}))$ is equal to
the Euler characteristic $\chi(mK_{X_t})$. Since the Euler
characteristic is a topological and thus a deformation invariant, it
follows that $h^0(X_t,\OOO_{X_t}(mK_{X_t}))$ is independent of $t$,
for $t$ small enough.
\end{proof}

This result is not general enough --- for a variety of general type
we know that $K_{X_0}$ is big, but not necessarily nef. To
generalize further, one is tempted to consider bundles that are big
and ``almost'' nef. However, if a bundle is not nef, it means that
there exists a curve $C\subset X$ such that $L\cdot C\le -1$. Thus
there are no integral divisors that are ``almost'' nef, and we are
led to consider divisors with rational coefficients.

\section{$\Q$-divisors and Kawamata-Viehweg's vanishing}
\begin{df}
A (Weil) $\Q$-divisor on $X$ is a formal linear combination $\sum
a_i D_i$, where $D_i$ are codimension one subvarieties of $X$, and
$a_i\in \Q$ are arbitrary coefficients. The rounding up and down of
a $\Q$-divisor are defined as $\lceil D\rceil:=\sum \lceil a_i\rceil
D_i$ and $\lfloor D\rfloor:=\sum \lfloor a_i\rfloor D_i$, by taking
respectively the least integer not smaller than, or the largest
integer not greater than $a_i$ for each $i$.
\end{df}

\begin{rem}
Notice that while all arithmetic operations with $\Q$-divisors are
functorial, i.e. commute with taking pullbacks and restrictions, the
rounding is not. Indeed, consider a line $\C$ in the plane $\C^2$
and the parabola $P$ touching it at the point $p$. Let $D=\frac 12
P$, so $\lfloor D\rfloor=0$. However, when we restrict $D$ to the
line and then round down, we get $\lfloor D|_{\C}\rfloor=\lfloor
2(p/2)\rfloor=p\ne 0=\lfloor D\rfloor |_\C$.
\end{rem}

Introducing $\Q$-divisors allows one to talk about $\Q$-divisors
that are close to an ample integral divisor. Moreover, one can
multiply $\Q$-divisors by a large number to clear all the
denominators, prove something for the resulting integral divisor,
and then take the appropriate root to recover properties of the
original divisor. However, it will turn out even more useful to work
with integral divisors that are close to an ample $\Q$-divisor.
Doing this allows one to perform some asymptotic constructions for
integral divisors, essentially by proving that {\it any} large
enough multiple $mL$ of a given integral divisor is uniformly
$\e$-close to {\it some} ample $\Q$-divisor, and then arguing that
then the integral divisor behaves almost as if it were ample, since
the error is essentially $\e/m$. This idea is indeed used in Siu's
proof of plurigenera, and setting it up rigorously is one of our
main goals. First, we will need to understand $\Q$-divisors better.

\begin{df}
Similarly to integral divisors, the support of a $\Q$-divisor $\sum
a_i D_i$ is the integral divisor $\sum_{\lbrace i| a_i\ne 0\rbrace}
D_i$. A $\Q$-divisor is said to have simple normal crossings
(s.n.c.) if its support has simple normal crossings, i.e. if all
$D_i$ are smooth, and whenever a number of $D_i$ intersect, their
normal vectors are linearly independent.
\end{df}

Simple normal crossings is the mildest singularity a divisor may
have. The pullback of a s.n.c. divisor under the blowup at a point
is also a s.n.c. divisor, and moreover restricting a divisor $D$ to
some $E$ such that the union of $D$ and $E$ is s.n.c. commutes with
rounding. The appropriate vanishing theorem for $\Q$-divisors is

\begin{thm}[Kawamata-Viehweg's vanishing, \cite{kaw},\cite{vie}]
Suppose $L$ is an integral divisor, numerically equivalent to (i.e.
its intersections with all effective curves are the same as those
of) $B+D$, where $B$ is a big and nef $\Q$-divisor, and $D$ is a
s.n.c. $\Q$-divisor such that $\lfloor D\rfloor=0$. Then the
vanishing holds for $L$, i.e. $H^i(X,\OOO_X(K_X+L))=0$ for $i>0$.
\end{thm}

This theorem makes precise what it means for an integral divisor $L$
to be ``close'' to a big and nef $\Q$-divisor $B$. We will not prove
this powerful result, but will rather give a different version of it
that will turn out to be equally useful.

\begin{cor} Suppose $L$ is an integral divisor, and $D$ is a s.n.c.
$\Q$-divisor such that $L-D$ is big and nef. Then
$H^i(X,\OOO_X(K_X+L-\lfloor D\rfloor))=0$ for $i>0$.
\end{cor}

Notice that restricting the Kawamata-Viehweg vanishing theorem to
the case of $D$ integral yields precisely the vanishing for big
and nef integral divisors. However, the improvement for
$\Q$-divisors is crucial for further applications. What happens is
essentially the following: the ample (or big and nef) cone is
indeed a cone, i.e. is invariant under scaling. Given a big
integral divisor that is not nef, it has to lie outside this cone,
and in fact has to be a fixed distance from it, as there is a
curve it intersects negatively. However, it may then happen that
its multiples will be sufficiently close to some ample
$\Q$-divisors, so that the Kawamata-Viehweg may be applied. Notice
that the set of divisors for which Kawamata-Viehweg's vanishing
holds is not invariant under scaling.

\section{Multiplier ideals and Nadel's vanishing}
In the previous section we have discussed the appropriate vanishing
condition for $\Q$-divisor with s.n.c. However, to be able to use
the vanishing in a wide variety of settings, we need to relax the
s.n.c. condition. The way to do this is to try to resolve a more
complicated singularity, and carrying this out naturally leads to
defining a multiplier ideal.

\begin{df} For a $\Q$-divisor $D$ on $X$, a log-resolution (also
called an embedded resolution) of the pair $(X,D)$ is a birational
map $\mu:X'\to X$ with $X'$ non-singular such that the divisor
$\mu^* D+{\rm except}(\mu)$ on $X'$ is s.n.c., where ${\rm
except}(\mu)$ denotes the exceptional locus of the birational map.
\end{df}

The log-resolutions are not unique, as we can blow up a s.n.c.
divisor again so that the preimage will still be s.n.c., or can do a
resolution in different ways. Hironaka proved that by doing an
appropriate sequence of blowups with smooth centers one can always
construct a log resolution. This is a hard result in resolution of
singularities, and as it falls far from the positivity questions
that we are interested in we don't discuss the proof here.

To get a version of the vanishing theorem that would work for
$\Q$-divisors with any kind of singularities, we would need to do a
log-resolution and then apply Kawamata-Viehweg for the resulting
s.n.c. divisor.

\begin{df} The relative canonical divisor of a rational map $\mu: X'\to
X$ is defined to be $K_{X'/X}:=K_{X'}-\mu^*K_X$. Though both $K_X$
and $K_{X'}$ are only defined as linear series, the relative
canonical is in fact an effective divisor: it is the locus where the
differential of the map $\mu$ is degenerate. The relative canonical
divisor is thus supported on ${\rm except}(\mu)$, and thus its
pushforward $\mu_*K_{X'/X}=\OOO_X$.
\end{df}

Suppose now that $\mu:X'\to X$ is a log-resolution of $(X,D)$. The
second form of Kawamata-Viehweg vanishing theorem on $X'$ then
gives for all $i>0$
$$
 0=H^i(X',\OOO_{X'}(K_{X'}+\lfloor\mu^*(L-D)\rfloor))
$$
$$
  =H^i(X',\OOO_{X'}(K_{X'/X}+\mu^*(K_X+L)-\lfloor \mu^* D\rfloor)).
$$
Now let us project this sheaf down to $X$ by $\mu$. If the higher
direct images $R^i$ vanish (and they in fact do, but we omit the
proof), by the projection formula we would also get the vanishing of
cohomology there:
$$
  0=H^i(X,\OOO_X(\mu_*(K_{X'/X}+\mu^*(K_X+L)-\lfloor\mu^*D\rfloor)))
$$
$$
  =H^i(X,\OOO_X(K_X+L)\otimes\mu_*(K_{X'/X}-\lfloor\mu^*D\rfloor)).
$$
This indicates what the appropriate correction for the vanishing
theorem should be.

\begin{df}[Esnault-Viehweg]
For a $\Q$-divisor $D$ on $X$ the multiplier ideal sheaf is defined
to be
$$
  \I(X,D):=\mu_*(K_{X'/X}-\lfloor\mu^*D\rfloor)
$$
Since $\mu_*(K_{X'/X})=\OOO_X$, the sheaf $\I(X,D)$ is an ideal
subsheaf of $\OOO_X$.
\end{df}

The discussion above indicates the direction one follows to prove
the appropriate vanishing theorem for $\Q$-divisors. The result is
as expected, but finishing the proof requires more work that we do
not show here --- see \cite{lazbook}, 9.4.B.

\begin{thm}[Nadel's vanishing] If $L$ is an integral
divisor and $D$ is a $\Q$-divisor such that $L-D$ is big and nef,
then
$$
  H^i(X,\OOO_X(K_X+L)\otimes\I(X,D))=0\quad{\rm for}\ i>0.
$$
\end{thm}

\begin{rem}

\noindent {\bf a)} From the definition it is not a priori clear that
$\I(X,D)$ does not depend on the choice of the log-resolution.
However, it can be proven that this is indeed the case.

\noindent {\bf b)} If $D$ is a $\Q$-divisor with s.n.c, the
log-resolution is the identity, the relative canonical bundle is
trivial, and thus $\I(X,D)=\OOO_X(-\lfloor D\rfloor)$

\noindent {\bf c)} If $D$ is an integral divisor, then there is no
rounding, and the definition commutes with pullback, so we get
$\I(X,D)=\OOO_X(-D)$.

\noindent {\bf d)} In view of the above, for the study of the
multiplier ideal all that matters is the fractional part of the
$\Q$-divisor: the integral part contributes a factor of minus
itself.
\end{rem}
\begin{df}
Since only the fractional part of the divisor is interesting for
multiplier ideal considerations, given any integral divisor $D$, it
is interesting to study the multiplier ideals $\I(X,cD)$ as $c$
increases from $0$ --- then the multiplier ideal is trivially
$\OOO_X$
--- to 1, when the multiplier ideal is $\OOO_X(-D)$ and thus non-trivial. The
log-canonical threshold of the pair $(X,D)$ is defined to be the
infinum of $c$ such that $\I(X,cD)$ is non-trivial. If this $c$ is
equal to one, i.e. if we have $\I(X,cD)=0$ for all $0<c<1$, then the
pair $(X,D)$ is called log-canonical.
\end{df}
\begin{ex}
To demonstrate how the multiplier ideals are computed, let us do
perhaps the simplest non-trivial example: $D$ is the sum of three
lines $\ell_1$, $\ell_2$ and $\ell_3$ intersecting at the point $p$
in the plane $\P^2$.

Let $\mu:X'\to\P^2$ be the blow up at $p$, and let $E$ be the
exceptional divisor on $X'$. The pullback of the divisor $cD$
$$
  \mu^*(cD)=c\mu^{-1}(\ell_1+\ell_2+\ell_3)+3cE
$$
is s.n.c.; thus $\mu$ is a log-resolution of $(\P^2,D)$. The
relative canonical divisor $K_{X'/\P^2}$ is equal to $E$. Thus for
$c<1/3$ we have $\lfloor \mu^* (cD)\rfloor=0$, while for $1/3\le
c<2/3$ we have $K_{X'/\P^2}-\lfloor \mu^* (cD)\rfloor=E-E=0$, so
that in both cases $\I(\P^2,cD)=\OOO_{\P^2}$. However, for $2/3\le
c<1$ we have
$$
\I(\P^2,cD)=\mu_*(K_{X'/\P^2}-\lfloor \mu^*
(cD)\rfloor)=\mu_*(-E)={\mathfrak m}_p
$$
is the maximal ideal of the point. Thus the log-canonical threshold
for $D$ is equal to $2/3$.

Doing the log-resolution explicitly (as three consecutive blowups),
one can show that for the case of $\ell=({\rm cusp}\ x^3=y^2)$ in
the plane we have $\I(\P^2,c\ell)=\OOO_{\P^2}$ for $c<5/6$, while
$\I(\P^2,c\ell)={\mathfrak m}_p$ for $5/6\le c$. It can be shown
that in general for the divisor of the curve $\lbrace
x^a=y^b\rbrace\subset \P^2$ the log-canonical threshold is equal to
$1/a+1/b$ --- this also works for three intersecting lines, which
are $\lbrace x^3=y^3\rbrace$.
\end{ex}

\section{Analytic approach to multiplier ideals and Nadel's vanishing}
In this section we will take a step back and indicate the analytic
construction of Nadel multiplier ideal sheaves in full generality,
and explain some of the analytic ideas that have been motivating
the work in the subject. This is a very rich analytic field, to
which we only present a simplified na\"\i ve introduction. The
reader is encouraged to consult \cite{siunote}, \cite{demsurv},
\cite{dembook} and references therein for a coherent rigorous
exposition.

In the analytic setting one starts with a line bundle $L$ on a
complex manifold $X$ of dimension $n$, and a Hermitian metric on it.
A Hermitian metric on a line bundle means simply a Hermitian scalar
product on all fibers. If we locally trivialize the bundle, by
choosing a basis vector $e(z)$ in each fiber (where
$z=(z_1,\ldots,z_n)$ is the local coordinate system on $X$), then
the Hermitian metric is determined by the value of the scalar
product $h(z):=(e(z),\overline e(z))$.

\begin{df}
The curvature of a Hermitian metric is the curvature of the unique
complex connection on the tangent bundle compatible with this
metric. It is a two-form on $X$ of type $(1,1)$, and is given in
coordinates as
$$
  \omega_h:=\frac{\sqrt{-1}}{2}\partial\overline\partial\log h:=
  \frac{\sqrt{-1}}{2}\sum_{i,j}\frac{\partial^2\log h(z)}{\partial
  z_i\partial\overline z_j}d z_i d\overline z_j.
$$
\end{df}
The curvature form is closed, and thus defines a class in
$H^{1,1}(X,\C)\cap H^2(X,\Z)$. By Poincar\'e duality (if we think of
$L$ as a linear combination of codimension one subvarieties) the
line bundle $L$ itself also defines a class in that cohomology
group, and the class of $\omega_h$ is linearly equivalent to $L$. In
the algebraic setting we are interested in computing the degrees of
the restriction of a line bundle to curves, or of its powers --- to
higher-dimensional subvarieties. Analytically the intersection
number with a subvariety $Y\subset X$ of dimension $m$ is
$$
  L^m\cdot Y=\int_Y(\omega_h)^m;
$$
notice that these integrals are independent of the choice of the
Hermitian metric $h$ on the line bundle $L$.

There is a natural invariant metric, the Fubini-Study metric, on the
tangent bundle to the projective space, which has constant positive
curvature. Thus if we use a very ample line bundle to embed some $X$
in a projective space, the restriction of the Fubini-Study metric to
the image will give a metric with positive curvature. This will also
hold for roots of very ample bundles, i.e. for ample bundles, and
thus the ampleness condition can be formulated analytically.

\begin{prop}
A line bundle $L$ is ample if and only if it admits such a Hermitian
metric $h$ that the corresponding curvature form is
positive-definite, i.e. $\omega_h(v,\overline v)>0$ for any non-zero
holomorphic tangent vector $v\in T^{1,0}(X)$.
\end{prop}

\begin{thm}[Kodaira's vanishing, analytic setting] If a line bundle
$L$ on $X$ admits a Hermitian metric of positive curvature, then the
cohomology groups $H^i(X,\OOO_X(K_X+L))$ are zero for all $i>0$.
\end{thm}

One is tempted to conjecture that a bundle is nef if and only if it
admits a Hermitian metric with non-negative curvature form. The
``if'' part --- if there is a Hermitian metric with non-negative
curvature, then the bundle is nef --- is indeed trivially true.
However, not all nef bundles do in fact admit metrics of
non-negative curvature. The correct criterion is
\begin{prop}
Fix a K\"ahler form $\nu$ on $X$ (i.e. a closed positive (1,1) form
on $X$, which does not need to be the curvature of some bundle);
then a line bundle $L$ on $X$ is nef if and only if for any $\e>0$
there exists a metric on $L$, whose curvature $\omega_\e$ is such
that $\omega_\e+\e\nu$ is strictly positive-definite.
\end{prop}

If a line bundle $L$ is not ample, then we cannot choose a metric on
it with positive-definite curvature form, so we would like to work
with a metric that is as close to being positive as possible. One
way to do this is to somehow bound from below the negativity of the
curvature form, but it turns out that the way that naturally leads
to multiplier ideals is rather to allow the metric to singularize,
while preserving the positivity.

Indeed, given a big line bundle (which is the case in our ultimate
goal here --- the invariance of plurigenera of general type, whence
the canonical bundle is big), we can use a high power of this bundle
defining a rational map to the projective space to pull back the
Fubini-Study metric. In doing this, we end up with a metric that is
well-defined, smooth, and has positive-definite curvature form away
from the indeterminacy locus of the rational map. On the
indeterminacy locus, however, the metric may singularize, i.e.
acquire a singularity as to become meromorphic there (in the
algebraic setting we don't encounter essential singularities). Thus
it will be natural in the following discussion to allow Hermitian
metrics with singularities.

\begin{rem}
It is amusing to note that when working with Riemann surfaces, one
also can either deal with a hyperbolic metric on a surface, i.e. a
metric of constant negative curvature, or can instead introduce a
flat Euclidean (zero-curvature) metric on a Riemann surface, which
would then be singular at a finite number of points. In a way what
we do now is a generalization and interpretation of this idea in
terms of metrics on general complex manifolds.
\end{rem}

Conceptually, we consider metrics that have zero curvature and
singularize along divisors. Suppose we are given an effective
irreducible $\Q$-divisor $D$, equal to $aH$ for some non-singular
subvariety $H\subset X$ of codimension one and some $a\in\Q_+$. Then
in local coordinate $z$ on $X$ near some point $p\in H$ the
subvariety $H$ is the zero locus of some function $f(z)$, the
vanishing order of which along $H$ is equal to one. Then locally we
can consider the singular Hermitian metric given by $|f(z)|^{-2a}$.
The curvature form of this metric,
\[
 \omega_f=\frac{\sqrt{-1}}{2}\partial\overline\partial (-2a\log |f|),
\]
is identically equal to zero away from $H$, and gives $a$ times
the delta-function of the subvariety $H$ (we recall that on $\C$
we have $\sqrt{-1}\partial\overline\partial\log |z|
=2\pi\delta(z)$).

A bit more precisely, what this means is the following. The
curvature ``form'' $\omega_f$ is in fact not a smooth form on $X$,
but rather has a singularity. However, given a smooth
$(n-1,n-1)$-form $\xi$ on $X$, the integral
$\int_X\omega_f\wedge\xi$ can still be computed. This means that
technically we should think of $\omega_f$ as a current of type (1,1)
--- an object dual to smooth $(n-1,n-1)$-forms on $X$. On the other
hand, $H$ can also be paired with $\xi$ by taking $\int_H\xi$. The
delta-function statement means simply that we have
$\int_X\omega_f\wedge\xi=a\int_H\xi$ for all $\xi$. Still a bit more
technically, we should note that this should only work locally, as
the defining function $f$ for $H$ can only be obtained locally, and
thus the above should rather be considered for smooth forms $\xi$ in
a small analytic neighborhood of a point $p\in H$. For the rigorous
discussion and details on currents we refer, for example, to
\cite{demsurv}.

The case of an irreducible divisor easily generalizes to the case of
a simple normal crossing divisor that is equal to $\sum a_i H_i$,
with $f_i$ being the (local) defining function for $H_i$. In this
case we consider the singular metric $\prod |f_i(z)|^{-2a_i}$.
Recall that for simple normal crossing divisors we know that the
multiplier ideal is $\I(X,D)=\OOO(-\lfloor D\rfloor)$. If we think
of this analytically, it means that the germ of the multiplier ideal
consists of locally defined holomorphic functions $F$ with vanishing
order along each of $H_i$ equal to at least $\lfloor a_i\rfloor$.
This is equivalent to $|F|^2/\prod |f_i(z)|^{2a_i}$ being locally
integrable near each of the $H_i$ and their intersections.

Analytically it is natural to attach the multiplier ideal to any
singular metric $e^{-2\phi}$, not only to those $\phi$ that come
from divisors. The technical condition is that the real function
$\phi$ must be plurisubharmonic, i.e. that whenever $\phi$ is finite
and smooth, its Laplacian $\sqrt{-1}\partial\overline\partial\phi$
is positive-definite, while $\phi$ is everywhere
upper-semi-continuous, and allowed to take the value of $-\infty$ at
some points.

\begin{df}
For a plurisubharmonic function $\phi$ on $X$ we define its
associated multiplier ideal $\I(\phi)\subset\OOO_X$ to be the sheaf
of germs of local holomorphic functions $F$ such that
$|F|^2e^{-2\phi}$ is locally integrable.
\end{df}

The discussion for effective normal crossing divisors can in fact be
applied to any divisor $D$, by passing to a small neighborhood where
all of its irreducible components become principal, and defining the
metric $e^{-2\phi_D}$ as above, irrespective of whether the
components are normal crossing or not.

\begin{prop}[see \cite{lazbook}, 9.3.D] In the algebraic setting,
the analytic and the algebraic multiplier ideals agree: for any
$\Q$-divisor $D$ we have $\I(X,D)=\I(\phi_D)$.
\end{prop}
We note, however, that the analytic multiplier ideal is defined in
more general setting, and allows one to potentially tackle the
non-algebraic situations as well.

Let us now discuss the analytic formulation and the intuition behind
Nadel's vanishing.
\begin{thm}[Analytic Nadel's vanishing, \cite{nadvani}]
Let $L$ be a line bundle on a projective algebraic variety $X$, with
singular metric $e^{-2\phi}$, such that its curvature current
$\omega_\phi$ dominates some smooth positive $(1,1)$-form $\nu$ on
$X$ as a current (i.e. for any positive $(n-1,n-1)$-form $\xi$ the
integral $\int_X (\omega_\phi-\nu)\wedge\xi>0$). Then the cohomology
groups $H^i(X,\OOO_X(K_X+L)\otimes\I(\phi))$ vanish for all $i>0$.
\end{thm}
\begin{proof}[Analytic idea of the proof, following \cite{siunote}]
The basic method is to try to smooth out the function $\phi$, while
controlling the smoothing in such a way that we can see that the
vanishing condition for the cohomology is preserved. The way to
smooth a function is to translate it by the flow of a holomorphic
vector field, and then average the translates (i.e. take the
integral of all translates for translation times from $-t$ to $t$
for some small $t$, and then divide the result by $2t$).

Let us embed $X$ in a large projective space $\P^N$ and then take a
generic projection $\pi:\P^N\to\P^n$. Let then
$U:=\pi^{-1}(\C^n)\subset X$, and let $Z:=X-U$ be the preimage
$\pi^{-1}(\P^{n-1})$. Let us choose a global trivialization for
$L|_U$.

The map $\pi$ restricted to $U$ is a branched cover. Let $Y$ be its
branch locus, so that the map $\pi:U-\pi^{-1}(Y)\rightarrow\C^n-Y$
is an unbranched cover. The vector fields $\frac{\partial}{\partial
z_j}$ on $\C^n-Y$ lift to holomorphic (this is why we needed to
eliminate the branching locus --- otherwise a branching singularity
could develop) vector fields on $U-\pi^{-1}(Y)$, which we denote by
$v_j$.

Notice that $Y$ is of codimension one in $\C^n$, and is the zero set
of some polynomial $F:\C^n\to\C$. Consider the family of compact
sets $K_a\subset \C^n$ for $a\ge 0$ such that $K_a$ is essentially
the set of points sufficiently far away from infinity and from $Y$,
and such that $\mathop{\cup}\limits_{a\ge 0} K_a=\C^n-Y$.
Technically, we can define the compact set $K_a:=\left\lbrace
z\Big|a\ge |F(z)|^{-2}+|z|^2\Big. \right\rbrace$, and denote its
preimage by $\Omega_a:=\pi^{-1}(K_a)$.

Now on $\Omega_a$ the vector fields $v_j$ are smooth and bounded (in
terms of $a$), and we can use them to smooth out $\phi$. In other
words, we can get {\it smooth} plurisubharmonic functions $\phi_\e$
on $\Omega_a$ monotonically decreasing to $\phi$ as $\e\to 0$. The
fact that the smoothings $\phi_\e$ are monotonically decreasing
follows from the sub-mean-value property, i.e. from the maximum
principle for plurisubharmonic functions.

Finally we are ready to think about the vanishing. We think of $H^i$
as the Dolbeault cohomology, i.e. an element of
$H^i(X,\OOO_X(K_X+L)\otimes\I(\phi))$ is a
$\overline\partial$-closed $(0,i)$ form $g$ on $X$ such that
$|g|^2e^{-2\phi}$ is integrable. To show that the cohomology is in
fact zero we need to show that there then exists some $(0,i-1)$-form
$u$ such that $\overline\partial u=g$. Then $g$ is exact and thus
represents the zero cohomology class.

Since we know that the smooth functions $\phi_\e\ge 0$ are
monotonically decreasing to $\phi$ and $|g|^2e^{-2\phi}$ is
integrable on $\Omega_a$, it means that $|g|^2e^{-2\phi_\e}$ is also
integrable on $\Omega_a$. Since everything is smooth and compact
now, we can find $u_{a,\e}$ such that $\overline\partial u_{a,\e}=g$
on $\Omega_a$ (the cohomology of $\C^n$ is trivial). Moreover, we
can choose $u_{a,\e}$ such that the integral $\int_{\Omega_a}
|u_{a,\e}|^2e^{-2\phi_\e}$, a.k.a. the $L^2$-norm of $u_{a,\e}$ on
$\Omega_a$ with weight $\phi_\e$, is bounded by some constant
independent of $a$ and $\e$.

The crucial step in the proof is the claim that this constant is
independent of $a$ and $\e$, and to prove this we use the fact that
$\omega_\phi$ dominates a smooth positive (1,1)-form $\nu$. Indeed,
using this fact we can bound both the growth of $u_{a,\e}$ near
$X-\Omega_a$, since there $\nu$ is also smooth, and the
approximation error in replacing $\phi$ by $\phi_\e$, as on compacts
by choosing $\e$ small enough we can ensure that $\phi_\e$ still
dominates $\nu$.

Thus finally we get a family of solutions $\overline\partial
u_{a,\e}=g$ on compacts $\Omega_a$ with bounded $L^2$ norms with
respect to weights $\phi_\e$. Taking first the limit as $\e\to 0$
and then the limit as $a\to\infty$, and using the fact that a
bounded normal family must converge, we finally show the existence
of the limit $u:=\lim\limits_{\e\to0,a\to\infty} u_{a,\e}$ with
bounded $L^2$ norm with respect to $\phi$, and it follows that
$\overline\partial u=g$ on $X$. We refer the reader to
\cite{nadvani} for the details of the original proof.
\end{proof}

\section{Multiplicity of multiplier ideals and Koll\'ar's theorem}
In this section we will study the multiplicity of divisors and the
non-triviality of the corresponding multiplier ideals, and will
prove Koll\'ar's theorem on the multiplicity of the theta divisor of
a principally polarized abelian variety.

\begin{prop}\label{mult}
Let $D$ be a $\Q$-divisor on $X$, let $n=\dim X$, and let $x\in X$
be a point. If the multiplicity ${\rm mult}_x D\ge n+p-1$, then
$\I(X,D)\subset{\mathfrak m}_x^p$.
\end{prop}
\begin{proof}
Let us construct a log-resolution $\mu:X'\to X$ of $(X,D)$ by first
blowing up the point $x$ and then doing whatever else is necessary.
Let us denote by $E$ the preimage in $X'$ of the exceptional divisor
of the first blow-up. We can compute ${\rm ord}_E(K_{X'/X}) =n-1$,
since $X'$ is computed by first blowing up at $x$, i.e. by inserting
a $\P^{n-1}$ instead of $x$.

The order to which the pullback of any divisor $F$ on $X$ contains
$E$ is equal to the multiplicity of vanishing of $F$ at $x$. Thus
${\rm ord}_E(\mu^*D)={\rm mult}_xD\ge n+p-1$, and in general
$\mu_*\OOO_{X'}(-aE)= {\mathfrak m}_x^a.$

Therefore
$$
  {\rm ord}_E(K_{X'/X}-\lfloor\mu^*D\rfloor)\le(n-1)-(n+p-1)=-p.
$$
It follows that
$$
  \OOO_{X'}(K_{X'/X}-\lfloor\mu^*D\rfloor) \subset \OOO_{X'}(-pE),
$$
and thus finally
$$
  \mu_*\OOO_{X'}(K_{X'/X}-\lfloor\mu^*D\rfloor)\subset\mu_*\OOO_{X'}
  (-pE)={\mathfrak m}_x^p.
$$
\end{proof}

This proposition shows that if a divisor has very high multiplicity
at some point, then the corresponding multiplier ideal sheaf is
non-trivial (i.e. not equal to $\OOO_X$). This can be generalized to
higher-dimensional subvarieties of $X$ --- the proof is analogous.

\begin{prop}
Suppose $Z\subset X$ is a subvariety of codimension $e$ such that
${\rm mult}_ZD\ge e+p-1$. Then $\I(X,D)\subset I_Z^{\langle
p\rangle}$, where $I_Z^{\langle p\rangle}$ is the symbolic power
--- the ideal of functions vanishing on $Z$ to order at least $p$.
\end{prop}

Now we will apply this to bound the multiplicity of theta divisors.

\begin{df} A principally polarized abelian variety $(A,\T)$ of dimension
$g$ is a complex torus $A$, i.e. a $g$-dimensional projective
variety with the structure of an abelian group on its points,
together with the choice of a principal polarization, i.e. an ample
line bundle $\T$ such that $h^0(A,\OOO_X(\T))=1$.
\end{df}

Principally polarized abelian varieties are a very classical object.
Classically they can be thought of as $A=\C^g/\Z^g+\tau\Z^g$, where
$\tau$ is a complex symmetric $g\times g$ matrix with
positive-definite imaginary part, and the theta function (the unique
up to a constant factor section of $\T$) is then given, for
$z\in\C^g$, by
$$
  \theta(z)=\sum\limits_{n\in\Z^g}\exp(\pi i(n,\tau n)+2\pi i(n,z)),
$$
where the transformation rule for $\theta(z)$ as we add to $z$ a
vector in the lattice $\Z^g+\tau\Z^g$ is what defines the bundle
$\T$ on $A$. Theta functions were studied extensively at least ever
since the works of Riemann. One natural question to ask is to
describe the vanishing locus of the theta function and its order of
vanishing. In particular one may wonder what is the maximal possible
multiplicity the theta function may have at a point. Despite this
being a question that already Riemann could ask, a classical
analytical solution is, to the best of our knowledge, still not
known. However, one can answer this question rather easily by using
multiplier ideals.

\begin{thm}[Koll\'ar, \cite{kolmult}]
The theta divisor (or theta function) cannot have multiplicity
greater than $g$ at any point of any principally polarized abelian
variety. More generally, for any $g$-dimensional principally
polarized abelian variety $A$ we have
$$
  {\rm dim} \lbrace x\in A|{\rm mult}_x\T\ge k\rbrace\le g-k.
$$
\end{thm}
\begin{proof}
We will use the above relation of the multiplicity of the divisor
and the non-triviality of the corresponding multiplier ideal sheaf.
In fact we will show below that the pair $(A,\T)$ is log-canonical.
Then the theorem would follow: indeed if at some point $x\in A$ we
had ${\rm mult}_x\T>g$, then for some $\e$ sufficiently close to 1
we would also have ${\rm mult}_x(\e\T)>g$ and thus by proposition
\ref{mult} the ideal $\I(\e\T)$ would be non-trivial, which would be
a contradiction. The bound for the dimension of the
high-multiplicity set is obtained analogously.

So let us prove that $\I(A,\e\T)=\OOO_A$ for all $0<\e<1$. Assume
the contrary: that for some $\e$ this is not the case, and then
denote by $Z$ the zero locus of $\I(A,\e\T)$. We must clearly have
then $Z\subset\T$. Consider the exact sequence
$$
  0\to\OOO_A(\T)\otimes\I(A,\e\Theta)\to\OOO_A(\Theta)\to\OOO_Z(\Theta)\to0,
$$
and look at the corresponding long exact sequence for cohomology.
Let us apply Nadel's vanishing to $A$ with the integral divisor
$L:=\Theta$ and the $\Q$-divisor $D:=\e\T$, so that indeed
$L-D=(1-\e)\T$ is ample. Since the canonical bundle $K_A$ is
trivial, we get for $i>0$
$$
  H^i(A,\OOO_A(K_A+L)\otimes\I(A,D))=H^i(A,\OOO_A(\T)\otimes\I(A,\e\T))=0.
$$

Thus from the long exact sequence we get the piece
$$
  H^0(A,\OOO_A(\T))\to H^0(Z,\OOO_Z(\T))\to 0,
$$
which must be exact. The first term of this sequence is just $\C$,
as $\T$ has a unique section, and thus the first map is zero, as we
are restricting the theta function to $Z$, which lies entirely in
its zero locus. Thus we must have $H^0(Z,\OOO_Z(\T))=0$; however,
this is impossible as we can always construct a section by
differentiating $\T$ in some direction sufficiently many times (and
using the fact that the derivative of a section $s$ of some line
bundle is the section of the same bundle when restricted to the zero
set of $s$). Alternatively we can prove it by translating the theta
divisor by some small $a\in A$, so that there would certainly be a
section (since $\T_a$ meets $\T$, and thus $Z\subset\T$, properly),
and then using semicontinuity to show that there is a section of
$\OOO_Z(\T)$ as well.

So we have arrived at a contradiction, and thus we must have $Z$
empty, so that $\I(A,\e\T)=\OOO_A$ for all $0<\e<1$, the pair
$(A,\T)$ is log-canonical, and thus the multiplicity of the theta
divisor is at most $g$ at all points.
\end{proof}

\section{Asymptotic methods and plurigenera}
In this and the next section we develop the necessary techniques and
explain the proof of Siu's theorem on invariance of plurigenera for
varieties of general type, following the algebraic exposition in
\cite{lazbook}, and referring to \cite{siuplur} for the original
analytic proof. Let us recall the statement:

\begin{thm}[Invariance of plurigenera]
If $X\to\Delta$ is a family of varieties of general type over the
unit disk $\Delta$, then for $t\in\Delta$ small enough the dimension
of $\Gamma(X_t,mK_{X_t})$ is independent of $t$, for all integer
$m$.
\end{thm}

As explained in proposition \ref{ampleplurigenus}, if we know that
the higher cohomologies of $mK_{X_0}$ vanish (and thus by
semicontinuity the higher cohomologies of $mK_{X_t}$ also vanish),
then the invariance of plurigenera follows from the invariance of
Euler characteristics. One way to try to ensure vanishing is by
applying Nadel's vanishing theorem. However, in Nadel's vanishing
the multiplier ideal enters, and thus we need to be able to control
it, which means controlling the singularities of divisors
$D_m\in|mK_{X_0}|$. If we try to do this directly for all $m$, the
log-resolutions get out of hand. This is when we are aided by the
fact that we are now working with $\Q$-divisors: philosophically
what we can try to do is resolve some $D_m$ for $m$ very large and
then take ``roots'' to resolve lower multiples of the canonical
bundle. To do this, let us introduce the asymptotic multiplier
ideals properly.

\begin{df}
The Iitaka dimension of a linear system $|L|$ on projective variety
$X$ of dimension $n$ is the number $\kappa(L):= \lim\limits_{
N\to\infty}\frac{\log h^0(X,\OOO_X(NL))} {\log N}$. Notice that by
definition $L$ is big if and only if $\kappa(L)=n$, and $X$ is of
general type if and only if $\kappa(K_X)=n$.
\end{df}

\begin{df}
For a linear system $|L|$ on $X$ with $\kappa(L)\ge 0$ the
asymptotic multiplier ideal is defined to be the direct limit
$$
  \I(X,c\,||L||):=\lim\limits_{N\to\infty}\I(X,\frac{c}{N}\,|NL|).
$$
Here the multiplier ideal of a linear system means that we choose a
general divisor $D_N\in|NL|$ and compute the multiplier ideal
$\I(X,\frac{c}{N}D_N)$, and thus the non-negativity of the Iitaka
dimension $\kappa(L)$ is needed to ensure the existence of such a
section $D_N$ for all $N$ sufficiently large. The existence of the
limit, and the fact that the ideals in the sequence grow as $m$
increases follows from the inclusion
$\I(X,\frac{c}{N}\,|NL|)\subset\I(X,\frac{c}{kN}\,|kNL|)$, which
holds for all integers $k>0$.
\end{df}

In general the asymptotic multiplier ideals differ from the usual
multiplier ideals. However, the following proposition provides us
with an example when they coincide.
\begin{prop}
If the ring of sections $\mathop{\oplus}\limits_{ m=1}^\infty
\Gamma(X,\OOO_X(mL))$ is finitely generated, then for $k\gg 0$ the
equality $\I(X, |mkL|)=\I(X,||mkL||)$ holds for all $m\ge 1$.
\end{prop}
\begin{proof}
Indeed, let us choose once and for all a log-resolution for all the
generators of the ring of sections simultaneously, and use it
independent of $m$ --- everything will get resolved, and thus we are
done.
\end{proof}

\begin{rem}
In view of this proposition, the asymptotic multiplier ideals for
linear systems with finitely generated rings of sections are easier
to deal with. However, for the case of the canonical linear system
that we are primarily concerned with, the finite generation is a
very hard open problem. Indeed, if finite generation were known, the
projectivization of the ring of sections --- the pluricanonical ring
--- would provide a canonical model for $X$ and thus prove
the minimal model program. Thus instead of using finite generation
to understand the multiplier ideals, one may instead try to use
multiplier ideals to tackle finite generation. This has not been
achieved yet, but the proof of the invariance of plurigenera can be
viewed as an indication that there may be hope in this approach.
\end{rem}

One way to understand a linear series is via its base locus. We have
the following easy observation

\begin{prop}\label{base}
The base ideal (the ideal dual to the base locus) of a linear system
is contained in the asymptotic multiplier ideal, i.e. ${\mathfrak
b}(|L|)\subset\I(X,||L||)$.
\end{prop}
\begin{proof}
The proposition is proven by fixing $k\gg 0$ large enough that it
computes the asymptotic multiplier ideal, i.e. such that
$\I(X,||L||)=\I(X,\frac{1}{k}|kL|)$), considering a log-resolution
$\mu:X'\to X$ resolving both $|L|$ and $|kL|$, noting that
${\mathfrak b}(|L|)^k\subset{\mathfrak b}(|kL|)$ and remembering
that the relative canonical class in the definition of the
multiplier ideal is an effective divisor.
\end{proof}

We now formulate the appropriate vanishing theorem for asymptotic
multiplier ideals.

\begin{thm}[see \cite{lazbook}, 11.2.12]\label{vanas}
Let $|L|$ be a linear system on $X$ with $\kappa(L)\ge 0$. Then

(a) For any big and nef integral divisor $A$ we have
$$
  H^i(X,\OOO_X(K_X+mL+A))\otimes\I(X,||mL||)=0
$$
for all $m\ge 1$ and all $i>0$.

(b) If we additionally assume $\kappa(L)=n$, i.e. that $L$ is big,
then it is enough to take $A$ nef, and not necessarily big, in the
above. In particular we can take $A=0$ to get in this case
$$
  H^i(X,\OOO_X(K_X+mL))\otimes\I(X,||mL||)=0.
$$
\end{thm}
Part (a) of this vanishing theorem follows from the fact that the
usual multiplier ideals stabilize to the asymptotic ideal, and then
applying the usual Nadel's vanishing. To prove part (b) more work is
needed, and it is achieved by essentially adapting the proof of
Nadel's vanishing to this case.

\begin{cor}\label{globgen}
In the setting of the theorem above, if $B$ is any ample globally
generated bundle, then
$$
  \OOO_X(K_X+mL+A+nB)\otimes\I(X,m||L||)
$$
is globally generated for all $m\ge 1$. If we additionally assume
$\kappa(L)=n$, then $A$ does not need to be big, and in particular
$A=0$ can be taken.
\end{cor}
\begin{proof}
From Castelnuovo-Mumford regularity (see \cite{lazbook}, 1.8) we
know that if for some very ample divisor $B$ and some sheaf $F$ we
have $H^i(X,F\otimes\OOO_X(-iB))=0$, then $F$ is globally generated.
Thus the global generation follows from the above vanishing theorem
for the asymptotic multiplier ideals.
\end{proof}

This in particular implies an interesting characterization of big
and nef bundles among all big bundles.
\begin{prop}
A big linear system $|L|$ is nef if and only if we have
$\I(X,||mL||)=\OOO_X$ for all $m\ge 1$.
\end{prop}
\begin{proof}[Proof of the ``if'' part]
Assume that the asymptotic multiplier ideals are trivial. Consider
then some very ample bundle $B$, and let $A=B$. Then the assumptions
of the above corollary are satisfied, and so
$$
\OOO_X(K_X+(n+1)B+mL)\otimes\I(m||L||)=\OOO_X(K_X+(n+1)B+mL)
$$
is big globally generated, and therefore nef. So for any effective
curve $C\subset X$ we have $(K_X+(n+1)B)\cdot C+mL\cdot C\ge 0$. By
taking the limit as $m\to\infty$ it then follows that $L\cdot C\ge
0$, and thus $L$ itself is nef.
\end{proof}

The theorems above are refinements of Nadel's vanishing, and thus
one can try to apply them to obtain deformation invariance of
plurigenera. The result one gets is the following

\begin{thm}[Siu's generalization of Levine's theorem]
If for some $k$ there exists a divisor $D\in|kK_{X_0}|$ such that
the pair $(X_0,D)$ is log-canonical, then the invariance of
plurigenera holds for all the plurigenera
$H^i(X_t,\OOO_{X_t}((m+1)K_{X_t}))$ for $m<k$.
\end{thm}
\begin{proof}
Indeed, $X_0$ is of general type, so $K_{X_0}$ is big, and we can
apply to it part (b) of theorem \ref{vanas} to get
$$
  H^i(X_0,\OOO_{X_0}((m+1)K_{X_0}))\otimes\I(X_0,||mK_{X_0}||)=0.
$$
By the log-canonicity assumption
$\I(X_0,\frac{m}{k}D)=\I(X_0,|mK_{X_0}|)=\OOO_{X_0}$ for any $m<k$,
and thus
$$
  \I(X_0,m||K_{X_0}||)=\OOO_{X_0}=\I(X_0,||mK_{X_0}||).
$$
Therefore the vanishing above becomes
$H^i(X_0,\OOO_{X_0}((m+1)K_{X_0}))=0$, and the invariance of
plurigenera $H^i(\OOO_{X_t}((m+1)K_{X_t}))$ for $m<k$ follows as
before in section \ref{ampleplurigenus} for the big and nef case.
\end{proof}

\section{Proof of Siu's theorem on the invariance of plurigenera}
We will now explain the proof of the invariance of plurigenera for
arbitrary varieties of general type. This is more delicate and more
technical than the considerations from the previous section, as the
vanishing by itself does not suffice. The trick will be to compare
the asymptotic multiplier ideals for multiples of $K_{X_0}$ and for
the multiples of $K_X$, then restricted to $X_0$. Indeed, given any
line bundle $L$ on the total space $X$ of the family, we can
consider the restriction maps
$$
  \phi_k:\Gamma(X,\OOO_X(kL))\to\Gamma(X_0,\OOO_{X_0}(kL_0))
$$
The images of these maps define a graded family of linear systems,
and thus we can consider the associated asymptotic multiplier ideal
as before, which we denote by
$$
  \I(X_0,||L||_0):=\lim\limits_{k\to\infty}\I\left(X_0,\frac{1}{k}
  \phi_{k}\Big(|\Gamma(X,\OOO_X(kL))|\Big)\right),
$$
where the subscript of $0$ after the notation of the asymptotic
multiplier ideal signifies the fact that we take sections of $L$
first and then restrict to the fiber over zero.

By working carefully with the definition of $\I(X_0,||L||_0)$, we
can see that (\cite{lazbook}, 11.5.5)
$$
  \Gamma\left(X_0,\OOO_{X_0}(K_{X_0}+L_0)\otimes\I(X_0,||L||_0)\right)\eqno(*)
$$
lies inside the image of the restriction map
$$
  \Gamma(X,\OOO_X(K_X+L))\to\Gamma(X_0,\OOO_{X_0}(K_{X_0}+L_0)),
$$
essentially since the asymptotic multiplier ideal comes from
restrictions of sections over $X$.

Thus for the case of $L=(m-1)K_X$ in the above, if we could show
that the whole space of sections of $K_{X_0}+L_0$ comes from the
sections in $(*)$, i.e. that
$$
  \Gamma(X_0,\OOO_{X_0}(mK_{X_0})\otimes\I(X_0,||(m-1)K_X||_0))=
  \Gamma(X_0,\OOO_{X_0}(mK_{X_0})),\eqno(**)
$$
the invariance of plurigenera would follow, as we would have all
sections of $mK_{X_0}$ lying in the image of restricting from
$mK_X$, i.e. we would then see that all sections of $mK_{X_0}$ can
be extended to sections of $mK_X$, and thus we would be done.

Thus what we need to prove is $(**)$, the fact that all the
sections of $mK_{X_0}$ vanish along $\I(X_0,||(m-1)K_X||_0))$. We
now notice that there is another asymptotic multiplier ideal that
can be naturally defined, which is just
$$
  \I(X_0,||mK_{X_0}||)=\lim\limits_{k\to\infty}\I\left(X_0,\frac{1}{k}
  \,\Big|\Gamma(X_0,\OOO_{X_0}(kmK_{X_0}))\Big|\right).
$$

By proposition \ref{base} (base locus is contained in the asymptotic
multiplier ideal) we know that the sections of $mK_{X_0}$ indeed do
vanish along this ideal, i.e. that
$$
  \Gamma\Big(X_0,\OOO_{X_0}(mK_{X_0})\otimes\I(X_0,||mK_{X_0}||)\Big)=
  \Gamma(X_0,\OOO_{X_0}(mK_{X_0})).
$$
Thus in particular if somehow we had
$$
 \I(X_0,||mK_{X_0}||)\subset\I(X_0,||(m-1)K_X||_0),
$$
the statement $(**)$ would follow and we would be done --- but of
course there is no reason for this inclusion to hold. However, there
is a weaker version that does hold, which still suffices to finish
the proof.

\begin{claim}\label{claim}
There exists an integer $a$ independent of $m$, and a pluricanonical
section $s$ of $aK_{X_0}$ with zero locus $D:=\lbrace
s=0\rbrace\subset X_0$, such that
$$
  \I(X_0,||mK_{X_0}||)(-D)\subset\I(X_0,||(m+a-1)K_X||_0).
$$
\end{claim}

The point of the claim is that though the two asymptotic multiplier
ideals we have defined are not the same, they are ``at most $a$ off
from each other'', and, very importantly, that this $a$ is
independent of the multiple $m$ of the canonical bundle that we
take. Given the claim, the proof of the invariance of plurigenera is
obtained in the following way.

Fix some $f\in\Gamma(X_0,mK_{X_0})$; we want to show that $f$
vanishes along the ideal $\I(X_0,||(m-1)K_X||_0)$. Consider the
section $f^N\cdot s\in\Gamma(X_0,(mN+a)K_{X_0})$ for $N$ very large.
It vanishes along $\I(X_0,||mNK_{X_0}||)$ (since $f^N$ vanishes
there), and also on $D$, since $s$ vanishes there. By the claim it
then follows that $f^N\cdot s$ vanishes along
$\I(X_0,||(mN+a-1)K_X||_0)$. We now use fact that for all multiplier
ideals we have $\I(X,||kL||)\subset\I(X,||L||)^k$ (this is known as
the subadditivity of the multiplier ideals, see \cite{lazbook}
9.5.B; notice that for base ideals, which are subideals of the
multiplier ideals by proposition \ref{base}, the inclusions
peculiarly goes the other way) this means that $f^N\cdot s$ vanishes
along $\I(X_0,||mK_X||_0)^N$. Since this is the case for all $N$, by
the claim this should imply that $f$ vanishes along
$\I(X_0,||(m-1)K_X||_0)$, once we choose $N\gg a$. Conceptually we
could actually try to say that $f$ vanishes on
$\I(X_0,||(m-\varepsilon)K_X||_0)$, but we do not need this.
Technically we need to use here the integral closure of ideals, and
the details can be found in \cite{lazbook}, 11.5.6. Thus the
invariance of the plurigenera for varieties of general type is
obtained once we prove the claim.

\begin{proof}[Proof of claim \ref{claim}]
The idea is that we want to add a very ample fixed piece to the
canonical divisor, so that everything becomes positive enough so
that we can handle the computation. In doing this we will be aided
by the so-called Kodaira's lemma: the statement that given any big
divisor $L$ and any divisor $F$ there is always a multiple of $L$
``bigger'' than $F$, i.e. that for some $a$ there exists an
effective divisor  in the linear system $|aL-F|$.

So choose a very ample bundle $B$ on $X$, positive enough so that
the linear system $F:=2K_X+(n+1)B$, where $n=\dim X_0$, is basepoint
free. Using Kodaira's lemma, choose $a$ such that there exists an
effective divisor $D$ in the linear series $|aK_X-F|$ --- this $D$
is then the zero set of some section $s\in\Gamma(X,\OOO_X(aK_X))$.
We will now show inductively that the statement of the claim holds
with these $a,s|_{X_0},$ and $D_0:=D|_{X_0}$ (in general the
subscript 0 will denote restricting divisors from $X$ to $X_0$).

To start the induction we need to show that
$$
  \I(X_0,||K_{X_0}||)(-D_0)\subset \I(X_0,||aK_X||_0).
$$
First notice that since the linear system $|2F|$ is basepoint free,
a general divisor $E\in|2F|$ is non-singular. Moreover, the
intersection of a general such $E$ with any subvariety is also
non-singular, and a generic $E$ intersects any smooth subvariety
transversely. Thus if we do a log-resolution for $(X,D)$, it also
gives a log-resolution for $(X,D+cE)$ for any constant $c$.
Therefore it follows from the definition of the multiplier ideal
that $\I(X,D)=\I(X,D+cE)$ for all $0<c<1$ --- this statement is
sometimes known as Koll\'ar-Bertini theorem. So finally the base of
induction is obtained by observing that
$$
  \I(X_0,||K_{X_0}||)(-D_0)\subset\OOO_{X_0}(-D_0)=\I(X_0,D_0)
$$
$$
  =\I(X_0,D_0+\frac{E_0}{2})=\I(X_0,|aK_{X_0}|)\subset\I(X_0,||aK_X||_0).
$$

Lastly, we need to prove the step of induction: that if the
inclusion of the claim holds for $k$, then it also holds for $k+1$.
Let us twist both sides of the claim by $\OOO_{X_0}((k+a)K_{X_0})$
--- it is enough to prove that the inclusion holds then. Writing
down what we get after this twist, on the left-hand-side of the
claim for $k+1$ we get
$$
  \OOO_{X_0}((k+a)K_{X_0}-D_0)\otimes\I(X_0,||(k+1)K_{X_0}||),\eqno
  (LHS)
$$
which is globally generated, since $aK_X-D$ is linearly equivalent
to $2K_X+(n+1)B$, and thus we can use corollary \ref{globgen} with
big $L:=K_X$, $m:=k+1$, and $A:=0$.

By the induction hypothesis we know that any section $u$ of
$(k+a)K_{X_0}$ that vanishes along $\I(X_0,||(k+1)K_{X_0}||)(-D_0)$
--- this is what the sections of (LHS) are --- must also vanish
along $\I(X_0,||(k+a-1)K_X||_0)$. Thus from the inclusion of the
space of sections $(*)$ in the image of the restriction map it
follows that this $u$ extends from the zero fiber, to give a section
$\tilde u$ of $(k+a)K_X$, which thus has to vanish along
$\I(X_0,||(k+a)K||_0)$. So we have shown that all the sections $u$
of (LHS) in fact lie in
$$
 \OOO_{X_0}(((k+a)K_{X_0})\otimes\I(X_0,||(k+a)K_X||_0).\eqno (RHS)
$$
But this is exactly the right-hand-side of the claim for $k+1$, and
since the sheaf (LHS) is globally generated, it means it is a
subsheaf of (RHS), which is precisely the statement of claim
\ref{claim} that we want, for $k+1$.
\end{proof}

\begin{proof}[Analytic viewpoint, following \cite{siuplur}]
The analytic version of the asymptotic multiplier ideal is obtained
as follows. For a bundle $L$ on $X$ with $\kappa(L)\ge 0$ one
chooses a basis $s_1^{(k)}\ldots s_{N_k}^{(k)}$ for
$\Gamma(X,\OOO_X(kL))$, takes the $k^{\text th}$ roots of these and
then arranges those, for all $k$, into a power series defining the
singular metric. Formally, one chooses a sequence $\lbrace
a_k\rbrace$ such that the sum
$$
  e^{-2\phi}:=\sum\limits_{k=1}^\infty a_k
  \sum\limits_{i=1}^{N_k}|s_i^{(k)}|^{2/k}
$$
converges uniformly. In the analytic setting the sequence $a_k$ can
be chosen arbitrarily, but should be the same for both
$\Gamma(X_0,\OOO_{X_0}(mK_{X_0}))$ and for
$\Gamma(X,\OOO_X(mK_X))|_{X_0}$, so that the analytic multiplier
ideals corresponding to the two resulting metrics (denote them by
$e^{-2\phi}$ and by $e^{-2\tilde\phi}$) that are constructed by
using these series can then be compared.

The analytic multiplier ideal is determined by the singular behavior
of the metric. Thus to compare two multiplier ideals one only needs
to compare the singularities of the metrics. Technically this means
that if one can show that the singularity of $e^{-2m\phi}$ is worse
than that of $e^{-2m\tilde\phi}$ by some fixed amount (i.e. the pole
order of the ratio is bounded, or something like that) for all $m$,
then the multiplier ideals corresponding to $e^{-2\phi}$ and
$e^{-2\tilde\phi}$ coincide. The technical result needed here to
show integrability is a generalization of an extension theorem of
Ohsawa-Takegoshi and Manivel.

The comparison of the singularities of $e^{-2m\phi}$ and of
$e^{-2m\tilde\phi}$, which is the analytic analog of claim
\ref{claim}, can be proven essentially the same way as in the
algebraic setting, using again the extension theorems, and also
Skoda's technique for generating elements of multiplier ideals.
\end{proof}
\begin{rem}
Siu later established in \cite{siuall} the technical analytical
result that was the sticking point in extending the analytic
techniques of \cite{siuplur} to arbitrary varieties, and thus proved
the deformation invariance of plurigenera for all varieties, not
necessarily of general type.
\end{rem}

\bibliographystyle{amsalpha}

\end{document}